\documentclass[oneside, English]{amsart}
\usepackage[T1]{fontenc}
\usepackage[latin1]{inputenc}
\usepackage{amssymb}

\makeatletter
  \theoremstyle{remark}
  \newtheorem{rem}{Remark}
 \theoremstyle{definition}
 \newtheorem{defn}{Definition}

\usepackage{latexsym}
\usepackage{amssymb,amsmath,amsxtra}
\newtheorem{theorem}{Theorem}[section]
\newtheorem{definition}{Definition}[section]

\newtheorem{remark}{Remark}[section]

\makeatother
\begin{document}

\title{An Integral Equation for Feynman's Operational Calcului}

\author{Lance Nielsen}

\curraddr{Department of Mathematics, Creighton University, Omaha, NE  68178}

\email{lnielsen@creighton.edu}
\begin{abstract}
In this paper we develop an integral equation satisfied by Feynman's operational calculi in formalism  of B. Jefferies and G. W. Johnson. In particular a ``reduced" disentangling is derived and an evolution equation of DeFacio, Johnson, and Lapidus is used to obtain the integral equation. After the integral equation is presented, we show that solutions to the heat and Schrodinger's equation can be obtained from the reduced disentangling and its integral equation. We also make connections between the Jefferies and Johnson development of the operational calculi and the analytic Feynman integral. 
\end{abstract}
\thanks{A substantial part of the research for this article was supported by the KRA grant (NSF 0354281)   at the University of Nebraska-Lincoln in May, 2006.}

\maketitle

\section{Introduction}

Feynman's operational calculus, originating with publication of the
paper \cite{Feyn51}, concerns itself with the formation of functions
of noncommuting operators. Indeed, even with functions as simple as
$f(x,y)=xy$, there is an ambiguity present in evaluating $f(A,B)$
if $A$ and $B$ do not commute. One is then left with the problem
of deciding, typically with a particular problem in mind, how best
to form a given function of noncommuting operators. One method of
dealing with this problem is to use the so-called Feynman indices.
That is, given operators $A$ and $B$, assign them indices. Then,
no matter how the product of $A$ and $B$ is written, the operator
with the smaller index acts first (to the right) of the operator with
the larger index. Hence, for example, given the product $\overset{1}{A}\overset{4}{B}$,
we can rewrite this as $\overset{4}{B}\overset{1}{A}$ since the index
of $A$ (= 1) is smaller than the index of $B$ (=4). Indeed, this
is an approach to the operational calculus taken by Maslov in \cite{Maslov}
and by Nasikinski, Shatalov, and Sternin in \cite{NSS}. 

In this paper we will follow the approach to the operational calculus
originated by Jefferies and Johnson in the series of papers \cite{JJ1,JJ2,JJ3,JJ4}
and expanded on in \cite{JJN}, \cite{JN1}, \cite{JNContDisc} and
others. In this approach to the calculus, the order of operators in
products is determined by the use of measures on intervals $[0,T]$.
One can use continuous measures as was done originally by Jefferies
and Johnson or one can use a mixture of continuous and discrete measures, see \cite{JNContDisc}.
We will use continuous measures in what follows. 

All of this being said, the question remains as to how one can use
measures to determine the order of operators in products. We will
start with a statement of Feynman's heuristic rules for the formation
of functions of noncommuting operators.

\vspace{8pt}

\noindent(1) Attach time indices to the operators to specify the
order of operators in products.

\noindent(2) With time indices attached, form functions of these
operators by treating them as though they were commuting.

\noindent(3) Finally, {}``disentangle'' the resulting expressions;
that is, restore the conventional ordering of the operators. 

\vspace{8pt}

As is well known, the central problem of the operational calculus
is the disentangling process. In his 1951 paper,  \cite{Feyn51}, Feynman points out
that ``The process is not always easy to perform and, in fact,
is the central problem of this operator calculus.'' 

We first address rule (1) above. This is, in fact, where we will see
measures coming into play. We first note that the operators
involved may come with time indices naturally attached. This is the case,
for example, with operators of multiplication by time-dependent potentials,
and also in connection with the Heisenberg representation in quantum
mechanics. However, it is also commonly the case that the operators 
used are independent of time. Given such an operator $A$, one
can (as Feynman most often did) attach time indices according to Lebesgue
measure as follows:
\[
A=\frac{1}{t}\int_{0}^{t}A(s)\, ds
\]
where $A(s):=A$ for $0\leq s\leq t$. While this device appears somewhat
artificial, it is extremely useful in many situations. Also,
it is worth noting that mathematical or physical considerations may
dictate that one use a measure different from Lebesgue measure in
a given situation. For example, if $\mu$ is a probability measure
on the interval $[0,T]$, and if $A$ is a linear operator, we may
write
\[
A=\int_{[0,T]}A(s)\,\mu(ds)
\]
where again $A(s):=A$ for $0\leq s\leq T$. Writing $A$ in this
fashion allows us to use the time variable $s$ to keep track of when
the operator $A$ acts. Indeed, consider two operators $A$ and $B$
and the product $A(s)B(t)$ (time indices have been attached). If
$s<t$, we have $A(s)B(t)=BA$ and if $t<s$, $A(s)B(t)=AB$; the
operator with the smaller or earlier time index acts before or to
the right of an operator with a larger or later time index. (We stress
here that these equalities are heuristic in nature.) For a more detailed
discussion of using measures to attach time indices to operators,
see the book \cite{JLBook}.

\begin{remark}
It is worth mentioning that using Dirac point mass measures to atttach time indices to the operators involved amounts to the use of the aforementioned Feynman indices.
\end{remark}

Concerning rules (2) and (3), we remark that, once time indices are
attached (so that an order of operation is specified), one can calculate
functions of the noncommuting operators by treating them as if they
actually do commute. Of course, such calculations are heuristic in
nature but the idea is that with time indices attached, one carries
out the necessary calculations giving no thought to the operator ordering
problem; the time indices will enable us to restore the desired ordering
of the operators once the calculations are finished -- this is the disentangling process.

The main subject of this article is to use the operational calculus
as put down by Jefferies and Johnson to derive an integral equation
satisfied by fully disentangled (or time ordered) operator expressions,
at least in a particular setting. Of critical importance in our derivation
will be the evolution equation found by DeFacio, Johnson, and Lapidus
in \cite{DJL} (it can also be found in \cite{JLBook}, Chapter 19).
In its initial form, the equation we will obtain will be rather general.
We will find, on the way to this equation, that the disentangled operator
we require will be somewhat different from the disentangling obtained
by Jefferies and Johnson. However, the difference will be that the disentangled
expression which will enable us to use the previously mentioned evolution
equation will contain fewer terms than the ``standard'' disentangled
expression one gets from Jefferies' and Johnson's formalism.

Once we obtain our integral equation we will consider several (related)
examples showing that the equation can be used to obtain solutions
to the heat equation and Schrodinger's equation. We will also use
the integral equation to establish some connections between disentangled
operators obtained in the spirit of Jefferies and Johnson and the
analytic Feynman integral defined using the Wiener integral. Both time
independent and time dependent potentials will be considered as well.
(See Chapters 13 and 15 of \cite{JLBook} for a very clear and detailed
discussion of the analytic Feynman integral.)

\section{The disentangling map}

Before the derivation of the integral equation, we define
the disentangling map in the time dependent setting. (The definition
is essentially identical in the time independent setting.) In doing
this we follow the initial definitions set out in \cite{NThesis}
and \cite{JJN}. 

\begin{remark}
As the reader will note, the algebras of functions defined below are
referred to as Banach algebras. We will not prove this assertion here
but will instead refer the reader to the paper \cite{JJ1} where the
proof is carried out for the time independent setting. However, as
noted in \cite{NThesis} and \cite{JJN}, the proof of this fact for
the time dependent setting is the same.
\end{remark}
\begin{remark}
We will assume throughout this section that the Banach space $X$
is separable.
\end{remark}

\begin{defn}
\label{def:A_Alg_Defn}Fix $T>0$. For $i=1,\ldots,k$, let $A_{i}:\,[0,T]\to\mathcal{L}(X)$
be maps that are measurable in the sense that $A_{i}^{-1}(E)$ is
a Borel set in $[0,T]$ for any strong operator open set $E\subset\mathcal{L}(X)$.
To each $A_{i}$ we associate a finite nonnegative continuous Borel measure $\mu_{i}$
on $[0,T]$ and require that, for each $i$,
\begin{equation}
r_{i}=\int_{[0,T]}\| A_{i}(s)\|_{\mathcal{L}(X)}\mu_{i}(ds)<\infty.\label{eq:Weight_Defn}
\end{equation}
We define, as in \cite{JJ1}, \cite{JJN}, and \cite{NThesis}, the
commutative Banach algebra $\mathbb{A}_{T}(r_{1},\ldots,r_{k})$ of
functions $f$ of $k$ complex variables that are analytic on the
open polydisk $\left\{ (z_{1},\ldots,z_{n})\,:\,|z_{i}|<r_{i},\, i=1,\ldots,k\right\} $
and such that the power series for $f(z_1,\ldots,z_k)$ converges on $\{(z_1,\ldots,z_k)\,:\, |z_i|=r_i,\,i=1\ldots,k\}$. (We emphasize that the weights we
are using here depend on the operator-valued functions as well as
on $T$ and on the measures.) The norm for this Banach algebra is
defined to be 
\begin{equation}
\| f\|_{\mathbb{A}_{T}}=\sum_{m_{1},\ldots,m_{k}=0}^{\infty}|a_{m_{1},\ldots,m_{k}}|r_{1}^{m_{1}}\cdots r_{k}^{m_{k}}\label{eq:A_Norm}
\end{equation}
where we use the Taylor series for $f$: 
\begin{equation}
f(z_{1},\ldots,z_{k})=\sum_{m_{1},\ldots,m_{k}=0}^{\infty}a_{m_{1},\ldots,m_{k}}z_{1}^{m_{1}}\cdots z_{k}^{m_{k}}.
\end{equation}

\end{defn}
In the next definition we define the \emph{disentangling algebra}.

\begin{defn}
To the algebra $\mathbb{A}_{T}$ we associate, as in \cite{JJ1}, \cite{JJN},
and \cite{NThesis}, a \emph{disentangling algebra} by creating formal
commuting objects $(A_{i}(\cdot),\mu_{i})\sptilde$, $i=1,\ldots,k$.
(These objects play the role of the indeterminants $z_{1},\ldots,z_{k}.$)
We define the disentangling algebra $\mathbb{D}_{T}\left((A_1,\mu_1)\sptilde,\ldots,(A_k,\mu_k)\sptilde \right)$
to be the collection of functions of the new indeterminants with the
same properties as the elements of the algebra defined in Definition
\ref{def:A_Alg_Defn}. However, rather than using the notation $(A_j,\mu_j)\sptilde$
below, we will often abbreviate to $A_j(\cdot)\sptilde$, especially
when carrying out calculations. The norm for $\mathbb{D}_{T}$ is
the same as that defined in (\ref{eq:A_Norm}) for the Banach algebra
$\mathbb{A}_{T}$ though we will refer to it as $\|\cdot\|_{\mathbb{D}_{T}}$
if a distinction needs to be made. 
\end{defn}
It is not hard to show that $\mathbb{A}_{T}$ and $\mathbb{D}_{T}$
are commutative Banach algebras which are isomorphic to one another
(see Propositions 1.1 - 1.3 in \cite{JJ1}).

For each $t\in[0,T]$ we now turn to the definition of the disentangling map 

\begin{equation}\label{E:MapNotation}
\mathcal{T}_{\mu_1,\ldots,\mu_k}^{t}\, : \,
\mathbb{D}_T\left((A_1(\cdot),\mu_1)\sptilde,\ldots,
(A_k(\cdot),\mu_k)\sptilde\right)\rightarrow \mathcal{L}(X).
\end{equation}

\noindent This will be done exactly as in \cite{JJ1}, \cite{JJN}, and \cite{NThesis}. In
order to state the next definition, which gives the action
of the disentangling map on monomials, we must first introduce some notation. (This
notation is essentially the same as used in \cite{JJ1}, \cite{JJN},
and \cite{NThesis}.) For a nonnegative integer $n$ and a permutation
$\pi\in S_{n}$, the set of all permutations of $1,\ldots,n$, we
define subsets $\Delta_{n}^{t}(\pi)$ of $[0,t]^n$ by 
\begin{equation}
\Delta_{n}^{t}(\pi)=\left\{ (s_{1},\ldots,s_{n})\in[0,t]^{n}:\,0<s_{\pi(1)}<\cdots<s_{\pi(n)}<t\right\} .\label{eq:DeltaSetDefinition}
\end{equation}
We next define, for nonnegative integers $n_{1},\ldots,n_{k}$ and
a permutation $\pi \in S_{n}$ with $n:=n_{1}+\cdots+n_{k}$, 

\begin{equation}\label{E:TimeOrderEqn}
        \tilde{C}_{\pi(i)}(s_{\pi(i)})=
        \begin{cases}
                A_{1}(s_{\pi(i)})\sptilde,  &\text{if $\pi(i) \in 
\{1,\ldots,n_{1}\}$},\\
                A_{2}(s_{\pi(i)})\sptilde,  &\text{if $\pi(i) \in 
                \{n_{1},\ldots,n_{1}+n_{2}\}$},\\
                \vdots \\
                A_{k}(s_{\pi(i)})\sptilde,  &\text{if $\pi(i)\in 
                \{ n_{1}+\cdots+n_{k-1}+1,\ldots,n\}$}. 
        \end{cases}
\end{equation}

\begin{remark}\label{R:CaseOfProbMeasures}
If the measures $\mu_1,\ldots,\mu_k$ are all probability measures it is straight forward to show that 
\begin{equation}\label{E:Commutative_World_Calculation}
\begin{split}
\left(A_1(\cdot)\sptilde\right)^{n_1} & \cdots \left(A_k(\cdot)\sptilde\right)^{n_k}=\\
&\sum_{\pi \in S_n} \int_{\Delta^t_n(\pi)} \tilde{C}_{\pi(n)}(s_{\pi(n)})\cdots 
\tilde{C}_{\pi(1)}(s_{\pi(1)}) \left(\mu_1^{n_1}\times \cdots \times \mu_k^{n_k} \right)(ds_1,\ldots,ds_n).
\end{split}
\end{equation}
This result is essential for the definition of the disentangling map.
\end{remark}

Now, for every $t\in[0,T]$, we define the action of the disentangling
map on monomials. 

\begin{definition}\label{D:MonomialDisentanglingMapDefinition}
Let $P_{t}^{n_{1},\ldots,n_{k}}\left(A_{1}(\cdot)\sptilde,\ldots,A_{k}(\cdot)\sptilde 
\right)=(A_{1}(\cdot)\sptilde)^{n_{1}}\cdots 
        (A_{k}(\cdot)\sptilde)^{n_{k}}$. \newline 
We define the action of the disentangling map on this monomial by (see Remark \ref{R:CaseOfProbMeasures})
\begin{equation}\label{E:TimeDepDefn}
 \begin{split}
 &\mathcal{T}_{\mu_{1},\ldots,\mu_{k}}^{t}
P_{t}^{n_{1},\ldots,n_{k}}\left(A_{1}(\cdot)\sptilde,\ldots,A_{k}(\cdot)\sptilde 
\right)\\
&=\mathcal{T}_{\mu_{1},\ldots,\mu_{k}}^{t}
 \left( \left(A_{1}(\cdot)\sptilde \right)^{n_{1}}\cdots
\left(A_{k}(\cdot)\sptilde \right)^{n_{k}} \right)\\
&:=\sum_{\pi \in S_{n}}\int_{\Delta_{n}^{t}(\pi)}
C_{\pi(n)}(s_{\pi(n)})\cdots C_{\pi(1)}(s_{\pi(1)})\left(\mu_{1}^{n_{1}}\times \cdots 
\times \mu_{k}^{n_{k}}\right)(ds_{1},\ldots,ds_{n})
\end{split}
\end{equation}
where the notation is as defined in \eqref{E:TimeOrderEqn} except that 
here we omit the tildes and consequently obtain the appropriate operator-valued functions in 
place of the formal commuting objects. 

Finally, for $f \in \mathbb{D}_{T}\left( 
(A_{1}(\cdot),\mu_{1})\sptilde,\ldots,(A_{k}(\cdot),\mu_{k})\sptilde\right)$ 
written as 
\begin{equation}\label{E:FunctionExpand}
        \begin{split}
            f\left(A_{1}(\cdot)\sptilde,\ldots,A_{k}(\cdot)\sptilde\right)=
            \sum_{n_{1},\ldots,n_{k}=0}^{\infty}c_{n_{1},\ldots,n_{k}}
            (A_{1}(\cdot)\sptilde)^{n_{1}}\cdots 
            (A_{k}(\cdot)\sptilde)^{n_{k}}
   \end{split}
\end{equation}
we define the action of the disentangling map on $f$ by
\begin{equation}\label{E:FullDisentDefn}
        \begin{split}           
&\mathcal{T}^{t}_{\mu_{1},\ldots,\mu_{k}}f\left(A_{1}(\cdot)\sptilde,\ldots,            
A_{k}(\cdot)\sptilde\right)\\
                &=:f_{t;\mu_{1},\ldots,\mu_{k}}\left(A_{1}(\cdot)\sptilde,\ldots,
                A_{k}(\cdot)\sptilde\right)\\
                &=\sum_{n_{1},\ldots,n_{k}=0}^{\infty}c_{n_{1},\ldots,n_{k}}
                \mathcal{T}^{t}_{\mu_{1},\ldots,\mu_{k}}                
P_{t}^{n_{1},\ldots,n_{k}}\left(A_{1}(\cdot)\sptilde,\ldots,A_{k}(\cdot)\sptilde\right).        
\end{split}
\end{equation}
\end{definition}

\begin{remark}
As is shown in \cite{NThesis} and \cite{JJN}, the disentangling
map is a linear contraction from the disentangling algebra into the
noncommutative Banach algebra of bounded linear operators on the Banach
space $X$. This differs somewhat from the time independent setting
of \cite{JJ1} where the disentangling map (defined exactly as above)
is shown to be a norm one contraction. As remarked in \cite{NThesis}
and \cite{JJN}, it is the presence of time dependent $\mathcal{L}(X)$
- valued functions that causes the map to be a contraction not necessarily
of norm one. 
\end{remark}

Next we define, using the sets $\Delta_{n}^{t}(\pi)$ defined above,
a set of permutations that will be very useful below. Let $n_{1},\ldots,n_{k}$
be nonnegative integers and consider the sets $\Delta_{n_{1}}^{t},\ldots,\Delta_{n_{k}}^{t}$
(where $\Delta_{j}^{t}:=\Delta_{j}^{t}(id)$ with ``$id$'' being
the identity permutation). Write, for each $\Delta_{n_{j}}^{t}$,
\begin{equation}
\Delta_{n_{j}}^{t}=\left\{ (s_{j,1},\ldots,s_{j,n_{j}})\in[0,t]^{n_{j}}:\,0<s_{j,1}<\cdots<s_{j,n_{j}}<t\right\} .
\end{equation}
(We therefore have $k$ ``blocks'' of completely time--ordered $n_{j}$ - tuples, $j=1,\ldots,k$.)
The set of permutations that we want are those permutations $\sigma$
of $\{1,\ldots,n\}$ that preserve the ordering of each block $s_{j,1},\ldots,s_{j,n_{j}}$.
We will denote this set of permutations by $\mathcal{P}_{n_{1},\ldots,n_{k}}$.
To be more specific, we require each $\sigma\in\mathcal{P}_{n_{1},\ldots,n_{k}}$
to preserve the order of the integers $n_{1}+\cdots+n_{j-1}+1,\ldots,n_{1}+\cdots+n_{j-1}+n_{j}$
in the list $\sigma(1),\sigma(2),\ldots,\sigma(n)$ though these integers
\emph{do not} have to appear consecutively in this list. In other words, the permutation $\sigma$ preserves the ordering in each block of time indices while at the same time putting the union of all the blocks of indices in the correct time order. The cardinality
of this set of permutations is 
\[
card(\mathcal{P}_{n_{1},\ldots,n_{k}})=\frac{(n_{1}+\cdots+n_{k})!}{n_{1}!\cdots n_{k}!}.
\]
Given $\sigma\in\mathcal{P}_{n_{1},\ldots,n_{k}}$, we let
\begin{multline}
(\Delta_{n_{1}}^{t}\times\cdots\times\Delta_{n_{k}}^{t})(\sigma)\\
=\left\{ (t_{1},\ldots,t_{n})\in\Delta_{n_{1}}^{t}\times\cdots\times\Delta_{n_{k}}^{t}:0<t_{\sigma(1)}<\cdots<t_{\sigma(n)}<t\right\} \qquad\qquad\label{eq:Delta_Product_Set}
\end{multline}
and note that, up to a set of measure zero (since we're using continuous measures),

\begin{equation}
\Delta_{n_{1}}^{t}\times\cdots\times\Delta_{n_{k}}^{t}=\bigcup_{\sigma\in\mathcal{P}_{n_{1},\ldots,n_{k}}}(\Delta_{n_{1}}^{t}\times\cdots\times\Delta_{n_{k}}^{t})(\sigma).
\end{equation}
A proof of this equality can be found in \cite{JLBook}.

It is easy to show that the disentangling map can be written using
a sum over $\mathcal{P}_{n_{1},\ldots,n_{k}}$ in place of the sum
over $S_{n}$ (see Proposition 2.5 of \cite{JJ1}). Indeed, for $f\in\mathbb{D}_{T}$, we may write

\begin{equation}\label{E:Alternate_Sum_Disent_Map}
\begin{split}
&\mathcal{T}_{\mu_1,\ldots,\mu_k}^{t}f(A_{1}\sptilde,\ldots,A_{k}\sptilde)=\\
&\sum_{n_1,\ldots,n_k=0}^{\infty} c_{n_1,\ldots,n_k} n_{1}!\cdots n_{k}! \sum_{\sigma \in \mathcal{P}_{n_1,\ldots,n_k}} 
\int_{(\Delta_{n_1}^{t}\times \cdots \times \Delta_{n_k}^{t})(\sigma)} \cdot \\
&C_{\sigma(n)}(s_{\sigma(n)})\cdots C_{\sigma(1)}(s_{\sigma(1)})
(\mu_{1}^{n_1} \times \cdots \times \mu_{k}^{n_k})(ds_1,\ldots,ds_n).
\end{split}
\end{equation}

\noindent A more complete discussion of this set of permutations,
with examples, can be found in \cite{DJL} or in \cite{JLBook}.

\section{The Integral Equation}

Let $\mathcal{H}$ be a separable Hilbert space. Let $A_{i}:[0,T]\to\mathcal{L}(\mathcal{H})$,
$i=1,\ldots,n$, be measurable in the sense of Definition \ref{def:A_Alg_Defn}.
To each $A_{i}$ associate a continuous nonnegative Borel measure $\mu_{i}$ on
$[0,T]$ and assume that 
\begin{equation}
r_{i}:=\int_{[0,T]}\| A_{i}(s)\|\,\mu_{i}(ds)<\infty\label{eq:Weights}
\end{equation}
for each $i$. Assume as well that $A_{i}(s)A_{i}(t)=A_{i}(t)A_{i}(s)$
for each $i$.

\begin{remark}
We do \emph{not}, however, assume that $A_{i}(s)A_{j}(t)=A_{j}(t)A_{i}(s)$
if $i\neq j$. 
\end{remark}

Finally, assume that the linear operator $-\alpha:\mathcal{H}\to\mathcal{H}$ generates
a $(C_{0})$ contraction semigroup of linear operators on $\mathcal{H}$. 

Using the nonnegative real numbers $r_{1},\ldots,r_{n}$ we construct
the commutative Banach algebra $\mathbb{A}_{T}\left(r_{1},\ldots,r_{n}\right)$.
Let $g\in\mathbb{A}_{T}\left(r_{1},\ldots,r_{n}\right)$ and define
\begin{equation}
f\left(z_{0},z_{1},\ldots,z_{n}\right)=e^{z_{0}}g(z_{1},\ldots,z_{n}).\label{eq:f_definition}
\end{equation}
Write
 \begin{equation}
g\left(z_{1},\ldots,z_{n}\right)=\sum_{m_{1},\ldots,m_{n}=0}^{\infty}g_{m_{1},\ldots,m_{n}}z_{1}^{m_{1}}\cdots z_{n}^{m_{n}}.\label{eq:g-Taylor_Series}
\end{equation}
 Associating Lebsegue measure $l$ with $-\alpha$ we have, using Definition 4.1 of \cite{JJN},
 \begin{multline}
\mathcal{T}_{l;\mu_{1},\ldots,\mu_{n}}^{t}f\left(-\alpha \sptilde,\,{A_{1}(\cdot)}\sptilde,\ldots,\,{A_{n}(\cdot)}\sptilde \right):=f_{l;\mu_{1},\ldots,\mu_{n}}^{t}\left(-\alpha,\, A_{1}(\cdot),\ldots,\, A_{n}(\cdot)\right)=\\
\sum_{m_{1},\ldots,m_{n}=0}^{\infty}g_{m_{1},\ldots,m_{n}}\sum_{\pi\in S_{m}}\int_{\Delta_{m}^{t}(\pi)}e^{-(t-s_{\pi(m)})\alpha}C_{\pi(m)}(s_{\pi(m)})\cdots C_{\pi(1)}(s_{\pi(1)})e^{-s_{\pi(1)}\alpha}\\
\left(\mu_{1}^{m_{1}}\times\cdots\times\mu_{n}^{m_{n}}\right)\left(ds_{1},\ldots,\, ds_{m}\right)\label{eq:Basic_Disentangling}
\end{multline}
 as the disentangling out to time $t\in[0,\, T]$. We can write this
expression as (see above or Proposition 2.5 of  \cite{JJ1})
\begin{multline}
f_{l;\mu_{1},\ldots,\mu_{n}}^{t}\left(-\alpha,\, A_{1}(\cdot),\ldots,\, A_{n}(\cdot)\right)=\\
\sum_{m_{1},\ldots,m_{n}=0}^{\infty}g_{m_{1},\ldots,m_{n}}m_{1}!\cdots m_{n}!\sum_{\pi\in\mathcal{P}_{m_{1},\ldots,m_{n}}}\int_{(\Delta_{m_1}^{t}\times \cdots \times \Delta_{m_n}^t)(\pi)}e^{-(t-s_{\pi(m)})\alpha}\\
C_{\pi(m)}(s_{\pi(m)})\cdots
C_{\pi(1)}(s_{\pi(1)})e^{-s_{\pi(1)}\alpha}\left(\mu_{1}^{m_{1}}\times\cdots\times\mu_{n}^{m_{n}}\right)\left(ds_{1},\ldots,\, ds_{m}\right)\label{eq:Modified_Disentangling}
\end{multline}
In particular we may write
\begin{multline}
\exp_{l;\mu_{1},\ldots,\mu_{n}}^{t}\left(-t\alpha+\sum_{j=1}^{n}\int_{[0,t]}A_{j}(s)\,\mu_{j}(ds)\right)=\\
\sum_{m_{1},\ldots,m_{n}=0}^{\infty}\sum_{\pi\in\mathcal{P}_{m_{1},\ldots,m_{n}}}\int_{\left(\Delta_{m_{1}}^{t}\times\cdots\times\Delta_{m_{n}}^{t}\right)(\pi)}e^{-(t-s_{\pi(m)})\alpha}C_{\pi(m)}(s_{\pi(m)})\cdots \\ \qquad \qquad \qquad \qquad C_{\pi(1)}(s_{\pi(1)})e^{-s_{\pi(1)}\alpha}
\left(\mu_{1}^{m_{1}}\times\cdots\times\mu_{n}^{m_{n}}\right)\left(ds_{1},\ldots,\, ds_{m}\right)\qquad\qquad\qquad\qquad\qquad\qquad\qquad\label{eq:Exponential_Disentangling}
\end{multline}
for the disentangling of the exponential function defined by $e^{z_{0}}e^{z_{1}+\cdots+z_{n}}$,
that is we are taking our function $g$ to be $g(z_{1},\ldots,z_{n})=e^{z_{1}+\cdots+z_{n}}$.
The disentangling displayed above in equation \eqref{eq:Exponential_Disentangling}
is that derived in the paper \cite{DJL} and can also be found in
\cite{JLBook}. 

While the expression seen in equation \eqref{eq:Modified_Disentangling} is the disentangling series for the function $f(A_1,\ldots,A_n)$, we will  not find it useful in obtaining our integral equation.  Indeed, in order to use the evolution equation from \cite{DJL} (see \eqref{eq:Evolution_Eqn} below) , we need the presence of the exponential function \eqref{eq:Exponential_Disentangling} in our disentangling. To change the form of our disentangling so that we can use the evolution equation, we first re-express the the power series coefficients using the Cauchy Integral Theorem and then slightly rewrite the operator functions to accommodate the use of the integral theorem. Using the integral form of $k!$, $k\in \mathbb{N}$, we arrive at equation \eqref{eq:Integral_Rep}.  This is equivalent to the disentangling seen in equation \eqref{eq:Modified_Disentangling} but it now contains the exponential function to which we will apply the aforementioned evolution equation. It is essentially equation \eqref{eq:Integral_Rep} that we will work with.  

As mentioned above, the first step in rewriting the disentangling $f_{l;\mu_{1},\ldots,\mu_{n}}^{t}$ in terms of the exponential function begins with the Cauchy Integral Theorem.  Indeed, we can use the Cauchy Theorem for derivatives to write the
coefficients in the power series \eqref{eq:g-Taylor_Series} as
 \begin{equation}
g_{m_{1},\ldots,m_{n}}=\frac{1}{(2\pi i)^{n}}\int_{|\xi_{1}|=r_{1}}\cdots\int_{|\xi_{n}|=r_{n}}g(\xi_{1},\ldots,\xi_{n})\xi_{1}^{-m_{1}-1}\cdots\xi_{n}^{-m_{n}-1}d\xi_{1}\cdots d\xi_{n}.\label{eq:Taylor_Coeffs}
\end{equation}
Now, replace each $A_{j}$ with $\frac{u_{j}}{\xi_{j}}A_{j}(\cdot)$.
Then, for $u_{j}\in[0,\infty)$ and $\xi_{j}\in\mathbb{C}\backslash\{0\}$
we obtain, provided we are able to interchange the order of summation
and integration as needed,
 \begin{multline}
f_{l;\mu_{1},\ldots,\mu_{n}}^{t}(-\alpha,\, A_{1}(\cdot),\ldots,\, A_{n}(\cdot))=\\
\frac{1}{(2\pi i)^{n}}\int_{|\xi_{1}|=r_{1}}\cdots\int_{|\xi_{n}|=r_{n}}\int_{[0,\infty)^{n}}g(\xi_{1},\ldots,\xi_{n})\xi_{1}^{-1}\cdots\xi_{n}^{-1}e^{-u_{1}}\cdots e^{-u_{n}}\cdot\\
\exp_{l;\mu_{1},\ldots,\mu_{n}}^{t}\left(-t\alpha+\sum_{j=1}^{n}\frac{u_{j}}{\xi_{j}}\int_{[0,t]}A_{j}(s)\,\mu_{j}(ds)\right)\, du_{1}\cdots du_{n}d\xi_{1}\cdots d\xi_{n}\label{eq:Integral_Rep}
\end{multline}
Of course, we have used the fact that
 \[
k!=\int_{[0,\infty)}u^{k}e^{-u}du
\]
for $k\in\mathbb{N}$. Hence, for $m_{1},\ldots,\, m_{n}\in\mathbb{N}$,
\[
m_{1}!\cdots m_{n}!=\int_{[0,\infty)^{n}}u_{1}^{m_{1}}\cdots u_{n}^{m_{n}}e^{-u_{1}}\cdots e^{-u_{n}}du_{1}\cdots du_{n}.
\]
Also, the presence of the factor $\frac{u_{j}}{\xi_{j}}$
with the operator $A_{j}$ in the exponential supplies the factor $u_{j}^{m_{j}}$
we need to obtain the necessary factorials in \eqref{eq:Integral_Rep}.
The $\xi_{j}$ in the denominator enables us to obtain the factors
$\xi_{j}^{-m_{j}-1}$ that are needed to give the power series coefficients.

We now verify that the above-mentioned interchanges of integration
and summation are valid. We use the vector--valued version of the standard theorem from analysis that states that if a sequence $\{g_n\}$ of functions is such that $\sum_{n=1}^{\infty}\int_{\Omega}|g_n|d\mu \in L^1$, then $\sum_{n=1}^{\infty} g_n \in L^1$ (see, for example, Corollary 12.33 of \cite{HS}). First note that
\begin{multline}
\sum_{m_{1},\ldots,m_{n}=0}^{\infty}\sum_{\pi\in\mathcal{P}_{m_{1},\ldots,m_{n}}}\bigg\{ \frac{1}{(2\pi i)^{n}}\int_{|\xi_{1}|=r_{1}}\cdots\int_{|\xi_{n}|=r_{n}}g(\xi_{1},\ldots,\,\xi_{n})\xi_{1}^{-m_{1}-1}\cdots\xi_{n}^{-m_{n}-1}\\
u_{1}^{m_{1}}\cdots u_{n}^{m_{n}}e^{-u_{1}}\cdots e^{-u_{n}}du_{1}\cdots du_{n}d\xi_{1}\cdots d\xi_{n}\bigg\} \int_{(\Delta_{m_{1}}^{t}\times\cdots\times\Delta_{m_{n}}^{t})(\pi)}e^{-(t-s_{\pi(m)})\alpha}\\
C_{\pi(m)}(s_{\pi(m)})\cdots
C_{\pi(1)}(s_{\pi(1)})e^{-s_{\pi(1)}\alpha}\left(\mu_{1}^{m_{1}}\times\cdots\times\mu_{n}^{m_{n}}\right)(ds_{1},\ldots,\, ds_{m}),\label{eq:Rearrange_Disent}
\end{multline}
being the disentangling series for $f\left(A_{1}(\cdot),\ldots,\, A_{n}(\cdot)\right)$,
converges absolutely in $\mathcal{L}(\mathcal{H})$ (\cite{JLBook}
or \cite{DJL}). It is clear that the scalar functions being summed/integrated
are continuous. The finite sum of integrals 
\begin{multline}
L_{m_{1},\ldots,m_{n}}:=\sum_{\pi\in\mathcal{P}_{m_{1},\ldots,m_{n}}}\int_{(\Delta_{m_{1}}^{t}\times\cdots\times\Delta_{m_{n}}^{t})(\pi)}e^{-(t-s_{\pi(m)})\alpha}C_{\pi(m)}(s_{\pi(m)})\cdots\\
C_{\pi(1)}(s_{\pi(1)})e^{-s_{\pi(1)}\alpha}\left(\mu_{1}^{m_{1}}\times\cdots\times\mu_{n}^{m_{n}}\right)(ds_{1},\ldots,\, ds_{m})\label{eq:Sum_of_Operators}
\end{multline}
is an operator in $\mathcal{L}(\mathcal{H})$ for $m_{1},\ldots, m_{n}\in\mathbb{N}$.
For fixed nonnegative integers $m_{1},\ldots, $ $m_{n}$, we can write
\begin{multline}
\int_{|\xi_{1}|=r_{1}}\cdots\int_{|\xi_{n}|=r_{n}}\int_{[0,\infty)^{n}}g(\xi_{1},\ldots,\,\xi_{n})\xi_{1}^{-m_{1}-1}\cdots\xi_{n}^{-m_{n}-1}u_{1}^{m_{1}}\cdots u_{n}^{m_{n}}\cdot\\
e^{-u_{1}}\cdots e^{-u_{n}}du_{1}\cdots du_{n}d\xi_{1}\cdots d\xi_{n}L_{m_{1},\ldots,m_{n}}\\
=\int_{[0,\infty)^{n}}\int_{|\xi_{1}|=r_{1}}\cdots\int_{|\xi_{n}|=r_{n}}g(\xi_{1},\ldots,\,\xi_{n})\xi_{1}^{-m_{1}-1}\cdots\xi_{n}^{-m_{n}-1}u_{1}^{m_{1}}\cdots u_{n}^{m_{n}}\cdot\\
e^{-u_{1}}\cdots e^{-u_{n}}d\xi_{1}\cdots d\xi_{n}du_{1}\cdots du_{n}L_{m_{1},\ldots,m_{n}}\\
=\int_{[0,\infty)^{n}}\int_{|\xi_{1}|=r_{1}}\cdots\int_{|\xi_{n}|=r_{n}}g(\xi_{1},\ldots,\,\xi_{n})\xi_{1}^{-m_{1}-1}\cdots\xi_{n}^{-m_{n}-1}u_{1}^{m_{1}}\cdots u_{n}^{m_{n}}\cdot\\
e^{-u_{1}}\cdots e^{-u_{n}}L_{m_{1},\ldots,m_{n}}d\xi_{1}\cdots d\xi_{n}du_{1}\cdots du_{n}.\label{eq:Integral_Interchange_Calc}
\end{multline}
The first equality follows by virtue of the standard Fubini theorem. The second equality follows from the fact that for $f:[a,\, b]\to\mathbb{C}$
integrable with respect to Lebesgue measure and for $T\in\mathcal{L}(\mathcal{H})$,
$\phi\in\mathcal{H}$, we have
\[
\int_{a}^{b}f(x)\, dxT\phi=\int_{a}^{b}f(x)T\phi\, dx=\int_{a}^{b}T(f(x)\phi)\, dx=T\left(\int_{a}^{b}f(x)\phi\, dx\right).
\]
The second, third, and fourth integrals just above are interpreted
as Bochner integrals. Further, note that
\begin{multline}
\sum_{m_{1},\ldots,m_{n}=0}^{\infty}\int_{[0,\infty)^{n}}\int_{|\xi_{1}|=r_{1}}\cdots\int_{|\xi_{n}|=r_{n}}g(\xi_{1},\ldots,\,\xi_{n})\xi_{1}^{-m_{1}-1}\cdots\xi_{n}^{-m_{n}-1}u^{m_{1}}\cdots u_{n}^{m_{n}}\cdot\\
e^{-u_{1}}\cdots e^{-u_{n}}L_{m_{1},\ldots,m_{n}}d\xi_{1}\cdots d\xi_{n}du_{1}\cdots du_{n}\label{eq:Series_of_Integrals}
\end{multline}
is the disentangling series for $f$ and so converges in $\mathcal{L}(\mathcal{H})$
- norm. Hence it follows that we are able to interchange the sum and
the integrals above, as asserted.

We may therefore write the disentangling of $f$ as shown in equation
\eqref{eq:Integral_Rep}. This integral expression gives us the disentangled
expression for the operator $f_{l;\mu_{1},\ldots,\mu_{n}}^{t}(-\alpha,$
$\, A_{1}(\cdot),\ldots,\, A_{n}(\cdot))$ corresponding to the function
\[
f(z_{0},\, z_{1},\ldots,\, z_{n})=e^{z_{0}}g(z_{1},\ldots,\, z_{n}).
\]
However, the utility of the integral expression \eqref{eq:Integral_Rep} for the disentangling
is not, for our purposes here, to further develop the operational
calculi for functions of the form above. Instead we will modify the
disentangling slightly in order to develop an integral equation for the operational
calculi (in the modified form). The reason for this modification is, even though we've now written 
$f_{l;\mu_{1},\ldots,\mu_{n}}^{t}$ with the disentangled exponential function, the presence of the factorial expression $m_1!\cdots m_n!$. These factorials would hinder the development of the appropriate integral equation as can clearly be seen when carrying out the calculations in equation \eqref{eq:Basi_Calculation}. 
We therefore discard the factor $m_{1}!\cdots m_{n}!$ from our disentangling. This means,
of course, omitting the integral over $[0,\infty)^{n}$ along with
the associated portions of the integrand. Consequently our disentangling series is "reduced" in the sense that we are summing over the smaller set of permutations $\mathcal{P}_{m_1,\ldots,m_n}$. What will be referred to henceforth as the "reduced disentangling" is, then, 
\begin{multline}
f_{(R);l;\mu_{1},\ldots,\mu_{n}}^{t}(-\alpha,\, A_{1}(\cdot),\ldots,\, A_{n}(\cdot))=\frac{1}{(2\pi i)^{n}}\int_{|\xi_{1}|=r_{1}}\cdots\int_{|\xi_{n}|=r_{n}}g(\xi_{1},\ldots,\,\xi_{n})\cdot\\
\xi_{1}^{-1}\cdots\xi_{n}^{-1}\exp_{l;\mu_{1},\ldots,\mu_{n}}^{t}\left(-t\alpha+\sum_{j=1}^{n}\frac{1}{\xi_{j}}\int_{[0,\, t]}A_{j}(s)\,\mu_{j}(ds)\right)\, d\xi_{1}\cdots d\xi_{n}.\label{eq:Reduced_Disentangling}
\end{multline}

\begin{rem}
From the structure of this expression it should be clear that it is
a sum of time - ordered operator products and so it is, by definition,
a disentangled expression in the sense of Feynman's rules. Indeed,
this expression contains the fully disentangled exponential function. 
\end{rem}
The reason we wish to work with the reduced disentangling will become
clear. From \cite{DJL} or \cite{JLBook} we know that the disentangled
exponential function \begin{equation}
E_{l;\mu_{1},\ldots,\mu_{n}}^{t}:=\exp_{l;\mu_{1},\ldots,\mu_{n}}^{t}\left(-t\alpha+\sum_{j=1}^{n}\frac{1}{\xi_{j}}\int_{[0,\, t]}A_{j}(s)\,\mu_{j}(ds)\right)\label{eq:Exp_Notation}\end{equation}
satisfies the evolution equation
\begin{equation}
E_{l;\mu_{1},\ldots,\mu_{n}}^{t}=e^{-t\alpha}+\sum_{j=1}^{n}\int_{[0,\, t]}e^{-(t-s)\alpha}\frac{1}{\xi_{j}}A_{j}(s)E_{l;\mu_{1},\ldots,\mu_{n}}^{s}\mu_{j}(ds).\label{eq:Evolution_Eqn}
\end{equation}
If we use this evolution equation in our reduced disentangling $f_{(R);l;\mu_{1},\ldots,\mu_{n}}^{t}$,
we obtain
\begin{equation}
\begin{split}
&f_{(R);l,\mu_{1},\ldots,\mu_{n}}^{t}(-\alpha,\, A_{1}(\cdot),\ldots,\, A_{n}(\cdot))=\frac{1}{(2\pi i)^{n}}\int_{|\xi_{1}|=r_{1}}\cdots\int_{|\xi_{n}|=r_{n}}g(\xi_{1},\ldots,\,\xi_{n})\cdot\\
&\xi_{1}^{-1}\cdots\xi_{n}^{-1}\left\{ e^{-t\alpha}+\sum_{j=1}^{n}\int_{[0,\, t]}e^{-(t-s)\alpha}\frac{1}{\xi_{j}}A_{j}(s)E_{l;\mu_{1},\ldots,\mu_{n}}^{s}\mu_{j}(ds)\right\} \, d\xi_{1}\cdots d\xi_{n}\\
&=\left\{ \frac{1}{(2\pi i)^{n}}\int_{|\xi_{1}|=r_{1}}\cdots\int_{|\xi_{n}|=r_{n}}g(\xi_{1},\ldots,\,\xi_{n})\xi_{1}^{-1}\cdots\xi_{n}^{-1}d\xi_{1}\cdots d\xi_{n}\right\} e^{-t\alpha}
\end{split}
\end{equation}
\begin{equation}
\begin{split}
&+\sum_{j=1}^{n}\frac{1}{(2\pi i)^{n}}\int_{|\xi_{1}|=r_{1}}\cdots\int_{|\xi_{n}|=r_{n}}g(\xi_{1},\ldots,\,\xi_{n})\xi_{1}^{-1}\cdots\xi_{j-1}^{-1}\xi_{j}^{-2}\xi_{j+1}^{-1}\cdots\xi_{n}^{-1}\cdot \nonumber \\
&\int_{[0,\, t]}e^{-(t-s)\alpha}A_{j}(s)E_{l;\mu_{1},\ldots,\mu_{n}}^{s}\mu_{j}(ds)d\xi_{1}\cdots d\xi_{n}\\
&\stackrel{=}{(*)}g(0,\ldots,\,0)e^{-t\alpha}+\sum_{j=1}^{n}\int_{[0,\, t]}e^{-(t-s)\alpha}A_{j}(s)\bigg\{ \frac{1}{(2\pi i)^{n}}\int_{|\xi_{1}|=r_{1}}\cdots\int_{|\xi_{n}|=r_{n}}\cdot\\
&g(\xi_{1},\ldots,\,\xi_{n})\xi_{1}^{-1}\cdots\xi_{j-1}^{-1}\xi_{j}^{-2}\xi_{j+1}^{-1}\cdots\xi_{n}^{-1}E_{l;\mu_{1},\ldots,\mu_{n}}^{s}d\xi_{1}\cdots d\xi_{n}\bigg\} \mu_{j}(ds).\label{eq:Basi_Calculation}
\end{split}
\end{equation}
(The inequality (*) above follows from the standard Fubini theorem.) 

Now, the expression in the braces in the last two lines of equation
\eqref{eq:Basi_Calculation} has the form of a reduced disentangling.
In order to proceed we need to identify the function being disentangled.
The key to this identification is contained in the integrals
\[
\frac{1}{(2\pi i)^{n}}\int_{|\xi_{1}|=r_{1}}\cdots\int_{|\xi_{n}|=r_{n}}g(\xi_{1},\ldots,\,\xi_{n})\xi_{1}^{-1}\cdots\xi_{j-1}^{-1}\xi_{j}^{-2}\xi_{j+1}^{-1}\cdots\xi_{n}^{-1}d\xi_{1}\cdots d\xi_{n}
\]
for $j=1,\ldots,\, n.$ When $E_{l;\mu_{1},\ldots,\mu_{n}}^{s}$ is
expanded in its disentangling series we obtain, for each $j$, the
following integral which is then written as the corresponding derivative
evaluated at the origin in $\mathbb{C}^{n}$ using Cauchy's Integral
Formula for derivatives:
\begin{multline}
\frac{1}{(2\pi i)^{n}}\int_{|\xi_{1}|=r_{1}}\cdots\int_{|\xi_{n}|=r_{n}}g(\xi_{1},\ldots,\,\xi_{n})\xi_{1}^{-m_{1}-1}\cdots\xi_{j-1}^{-m_{j-1}-1}\xi_{j}^{-m_{j}-2}\cdot\\
\xi_{j+1}^{-m_{j+1}-1}\cdots\xi_{n}^{-m_{n}-1}d\xi_{1}\cdots d\xi_{n}\qquad\\
=\frac{\frac{\partial^{m_{1}+\cdots+m_{j-1}+(m_{j}+1)+m_{j+1}+\cdots+m_{n}}g}{\partial z_{1}^{m_{1}}\cdots\partial z_{j-1}^{m_{j-1}}\partial z_{j}^{m_{j}+1}\partial z_{j+1}^{m_{j+1}}\cdots\partial z_{n}^{m_{n}}}(0,\ldots,0)}{m_{1}!\cdots m_{j-1}!(m_{j}+1)!m_{j+1}!\cdots m_{n}!}\qquad\qquad\qquad\qquad\qquad\label{eq:g_Coefficients}
\end{multline}
In order to determine the scalar function giving the disentangling
with the Taylor coefficients just above, we consider the corresponding
power series for each $j=1,\ldots,n$: 
\begin{multline}
\sum_{m_{1},\ldots,m_{n}=0}^{\infty}\left(\frac{\frac{\partial^{m_{1}+\cdots+m_{j-1}+(m_{j}+1)+m_{j+1}+\cdots+m_{n}}g}{\partial z_{1}^{m_{1}}\cdots\partial z_{j-1}^{m_{j-1}}\partial z_{j}^{m_{j}+1}\partial z_{j+1}^{m_{j+1}}\cdots\partial z_{n}^{m_{n}}}(0,\ldots,0)}{m_{1}!\cdots m_{j-1}!(m_{j}+1)!m_{j+1}!\cdots m_{n}!}\right)\cdot \\
z_{1}^{m_{1}}\cdots z_{j-1}^{m_{j-1}}z_{j}^{m_{j}}z_{j+1}^{m_{j+1}}\cdots z_{n}^{m_{n}}
=\sum_{m_{1},\ldots,m_{j-1},m_{j+1},\ldots,m_{n}=0}^{\infty }\cdot \qquad \qquad\\
\left\{ \frac{1}{z_{j}}\sum_{m_{j}=0}^{\infty}\left(\frac{\frac{\partial^{m}g(0,\ldots,0)}{\partial z_{1}^{m_{1}}\cdots\partial z_{j-1}^{m_{j-1}}\partial z_{j}^{m_{j}}\partial z_{j+1}^{m_{j+1}}\cdots\partial z_{n}^{m_{n}}}}{m_{1}!\cdots m_{j-1}!m_{j}!m_{j+1}!\cdots m_{n}!}\right)z_{1}^{m_{1}}\cdots z_{j-1}^{m_{j-1}}z_{j}^{m_{j}}z_{j+1}^{m_{j+1}}\cdots z_{n}^{m_{n}}\right.\\
\left.-\frac{1}{z_{j}}\left(\frac{\frac{\partial^{m_{1}+\cdots+m_{j-1}+m_{j+1}+\cdots+m_{n}}g(0,\ldots0)}{\partial z_{1}^{m_{1}}\cdots\partial z_{j-1}^{m_{j-1}}\partial z_{j+1}^{m_{j+1}}\cdots\partial z_{n}^{m_{n}}}}{m_{1}!\cdots m_{j-1}!m_{j+1}!\cdots m_{n}!}\right)z_{1}^{m_{1}}\cdots z_{j-1}^{m_{j-1}}z_{j+1}^{m_{j+1}}\cdots z_{n}^{m_{n}}\right\} 
\end{multline}
\begin{multline}
=\frac{1}{z_{j}}\sum_{m_{1},\ldots,m_{n}=0}^{\infty}\left(\frac{\frac{\partial^{m}g(0,\ldots,0)}{\partial z_{1}^{m_{1}}\cdots\partial z_{n}^{m_{n}}}}{m_{1}!\cdots m_{n}!}\right)z_{1}^{m_{1}}\cdots z_{n}^{m_{n}}-\qquad\qquad\qquad\qquad\qquad\qquad\qquad\qquad \nonumber \\
\frac{1}{z_{j}}\sum_{m_{1},\ldots,m_{j-1},m_{j+1},\ldots,m_{n}=0}^{\infty}\left(\frac{\frac{\partial^{m_{1}+\cdots+m_{j-1}+m_{j+1}+\cdots+ m_{n}}g(0,\ldots,0)}{\partial z_{1}^{m_{1}}\cdots\partial z_{j-1}^{m_{j-1}}\partial z_{j+1}^{m_{j+1}}\cdots\partial z_{n}^{m_{n}}}}{m_{1}!\cdots m_{j-1}!m_{j+1}!\cdots m_{n}!}\right)
\cdot \\
z_{1}^{m_{1}}\cdots z_{j-1}^{m_{j-1}}z_{j+1}^{m_{j+1}}\cdots z_{n}^{m_{n}}
=\frac{1}{z_{j}}\left(g(z_{1},\ldots,\, z_{n})-g(z_{1},\ldots,z_{j-1},\,0,\, z_{j+1},\ldots,\, z_{n})\right)\qquad\qquad\qquad\qquad\qquad\qquad\qquad\label{eq:Int_Eqn_Func}
\end{multline}
With this done, we may now record our integral equation:
\begin{theorem}\label{T:IntegralEquation}
Under the hypotheses stated at the beginning of the section, the reduced disentangling 
$f^t_{(R);l;\mu_1,\ldots,\mu_n}$ satisfies the equation
 \begin{equation}
 \begin{split}
&f_{(R);l;\mu_{1},\ldots,\mu_{n}}^{t}(-\alpha,\, A_{1}(\cdot),\ldots,\, A_{n}(\cdot))
=g(0,\ldots,\,0)e^{-t\alpha}+\\
&\sum_{j=1}^{n}\int_{[0,\, t]}e^{-(t-s)\alpha}A_{j}(s)\Phi_{(R);l;\mu_{1},\ldots,\mu_{n}}^{j,\, s}(-\alpha,\, A_{1}(\cdot),\ldots,\, A_{n}(\cdot))\,\mu_{j}(ds)\label{eq:Integral_Equation}
\end{split}
\end{equation}
where
\begin{equation}
\begin{split}
\Phi^{j}(z_{0},\, z_{1},\ldots,\, z_{n})&=\frac{e^{z_{0}}}{z_{j}}\left(g(z_{1},\ldots,\, z_{n})-g(z_{1},\ldots,\, z_{j-1},0,\, z_{j+1},\ldots,\, z_{n})\right)\\
&=\frac{1}{z_j}\left(f(z_0,z_1,\ldots,z_n)-f(z_0,z_1,\ldots,z_{j-1},0,z_{j+1},\ldots,z_n)\right)
\end{split}
\label{eq:Phi_Definition}
\end{equation}
for $j=1,\ldots,n$.
\end{theorem}

\section{Examples}

In order to investigate the integral equation above, we consider some
examples. Before proceeding to our examples, however, we define the
mixed norm space $L_{\infty1;l}$. We say that $V\in L_{\infty1;\eta}([0,\, T]\times\mathbb{R}^{d})$
if \[
\| V\|_{\infty1;\eta}:=\int_{[0,\, T]}\| V(s,\cdot)\|_{\infty}\eta(ds)\]
where $\eta\in M([0,\, T])$. Given $V\in L_{\infty1;l}([0,\, T]\times\mathbb{R}^{d})$
($l$ is Lebesgue measure) it follows that $M_{V}$, the operator
of multiplication by $V$, is a bounded linear operator on $L^{2}(\mathbb{R}^{d})$
for Lebesgue almost every $s\in[0,\, T]$. The norm of $M_{V}$ is
$\| M_{V}\|=\| V(s,\cdot)\|_{\infty}$. (On can find a discussion of the mixed norm space in \cite{JLBook}.)

We now proceed with our first example. Let $A_{1}$ be the operator
of multiplication by $V:[0,\, T]\times\mathbb{R}^{d}\to\mathbb{R}$,
$V\in L_{\infty1;l}$. Take $\alpha=H_{0}=\frac{1}{2}\Delta$. Define\begin{equation}
f(z_{0},\, z_{1})=e^{z_{0}}\left(z_{1}g_{1}(z_{1})+g_{1}(0)\right)\label{eq:Two_Var_Fn}\end{equation}
for $g_{1}(z_{1})$ analytic on the disk $D(0,\| V\|_{\infty1;l})$
and continuous on its boundary. (Here we are replacing the function
$g(z)$ above with the function $z_{1}g_{1}(z_{1})+g_{1}(0)$.) It is well known
that, for $t>0$ and $\psi\in L^{2}(\mathbb{R}^{d})$, 
\begin{equation}
\left(e^{-tH_{0}}\psi\right)(x)=(2\pi t)^{-d/2}\int_{\mathbb{R}^{d}}\psi(u)\exp\left(\frac{-\| x-u\|^{2}}{2t}\right)\, du.\label{eq:Heat_Semi_gp}
\end{equation}
To proceed with our calculations we write the Taylor series for $g_1$:
\begin{equation}
g_{1}(z_{1})=\sum_{m=0}^{\infty}a_{m}z_{1}^{m}.\label{eq:g1_PWR_Series}
\end{equation}
To do the calculation on the right hand side of the integral equation
\eqref{eq:Integral_Equation}, we first note that here $n=1$ and associating
Lebesgue measure to $A_{1}$, the operator of multiplication by $V$,
we have
\begin{equation}
f_{(R);l;l}^{t}(-H_{0},V)=g_{1}(0)e^{-tH_{0}}+\int_{0}^{t}e^{-(t-s)H_{0}}V(s)\Phi_{(R);l;l}^{1,s}(-H_{0},\, V)\, ds\label{eq:One_Dim_IE}
\end{equation}
 for our integral equation. From the definition of $\Phi^{j}$ we have,
using $z_{1}g_{1}(z_{1})+g_{1}(0)$ in place of $g(z_{1})$, 
\begin{equation}
\Phi^{1}(z_{0},\, z_{1})=e^{z_{0}}\cdot\frac{1}{z_{1}}(z_{1}g_{1}(z_{1})+g_{1}(0)-g_{1}(0))=e^{z_{0}}g_{1}(z_{1}).\label{eq:New_g}
\end{equation}
Hence the reduced disentangling $\Phi_{(R);l;l}^{1,s}(-H_{0},\, V)$
is that of the function $e^{z_{0}}g_{1}(z_{1})$ for $z_{0}=-H_{0}$
and $z_{1}=V$. Using the power series for $g_{1}$ we have,
given $\psi\in L^{2}(\mathbb{R}^{d})$, 
\begin{equation}
\begin{alignedat}{1} & (\Phi_{(R);l;l}^{1,s}(-H_{0},V)\psi)(x)\\
= & \sum_{m=0}^{\infty}a_{m}\int_{\Delta_{m}(s)}\left(e^{-(s-t_{m})H_{0}}V(t_{m})\cdots V(t_{1})e^{-t_{1}H_{0}}\psi\right)(x)\, dt_{1}\cdots dt_{m}\\
= & \sum_{m=0}^{\infty}a_{m}\int_{\Delta_{m}(s)}\left((2\pi)^{m+1}(s-t_{m})\cdots(t_{2}-t_{1})t_{1}\right)^{-d/2}\int_{\mathbb{R}^{(m+1)d}}V(x_{m},t_{m})\cdot\\
 & V(x_{m-1},t_{m-1})\cdots V(x_{2},t_{2})V(x_{1},t_{1})\psi(x_{0})\exp\left[-\sum_{j=1}^{m+1}\frac{\| x_{j}-x_{j-1}\|^{2}}{2(t_{j}-t_{j-1})}\right]\cdot\\
 & \qquad\qquad dx_{0}dx_{1}\cdots dx_{m-1}dt_{1}dt_{2}\cdots dt_{m}\end{alignedat}
\label{eq:Initial_Phi_Calc}
\end{equation}
where $x_{m+1}=x$ and $t_{m+1}=s$ in the sum inside the exponential
above. We now calculate\[
e^{-(t-s)H_{0}}V(s)\Phi_{(R);l;l}^{1,s}(-H_{0},V)\psi.\]
Using continuity of the semigroup we may write\begin{equation}
\begin{alignedat}{1} & \left(e^{-(t-s)H_{0}}V(s)\Phi_{(R);l;l}^{1,s}(-H_{0},V)\psi\right)(x)\\
= & \sum_{m=0}^{\infty}a_{m}\int_{\Delta_{m}(s)}\left((2\pi)^{m+2}(t-s)(s-t_{m})(t_{m}-t_{m-1})\cdots(t_{2}-t_{1})t_{1}\right)^{-d/2}\cdot\\
 & \int_{\mathbb{R}^{(m+2)d}}V(x_{m+1},\, t_{m+1})V(x_{m},\, t_{m})\cdots V(x_{2},\, t_{2})V(x_{1},\, t_{1})\psi(x_{0})\cdot\\
 & \exp\left[-\sum_{j=1}^{m+2}\frac{\| x_{j}-x_{j-1}\|^{2}}{2(t_{j}-t_{j-1})}\right]\, dx_{0}dx_{1}\cdots dx_{m}dt_{1}dt_{2}\cdots dt_{m}\end{alignedat}
\label{eq:Further_Calcuation}\end{equation}
where $x_{m+2}=x$, $t_{m+2}=t$ and $t_{m+1}=s$. For any $y\in C_{0}^{t}$,
the classical Wiener space of paths, define $G:\, C_{0}^{t}\to\mathbb{R}$
by \begin{equation}
G(y):=V(y(t_{1}),t_{1})\cdots V(y(t_{m+1}),t_{m+1})\psi(y(t)).\label{eq:G_DEFN}\end{equation}
Then\begin{alignat}{1}
\int_{C_{0}^{t}}G(y+x)\, m(dy)=\qquad\qquad\qquad\qquad\qquad\qquad\qquad\qquad\nonumber \\
\left((2\pi)^{m+2}(t-s)(s-t_{m})\cdots(t_{2}-t_{1})t_{1}\right)^{-d/2}\int_{\mathbb{R}^{(m+2)d}}V(x_{m+1},s)V(x_{m},t_{m})\nonumber \\
\cdots V(x_{2},t_{2})V(x_{1},t_{1})\psi(x_{0})\exp\left[-\sum_{j=1}^{m+2}\frac{\| x_{j}-x_{j-1}\|^{2}}{2(t_{j}-t_{j-1})}\right]\, dx_{0}\cdots dx_{m}dt_{1}\cdots dt_{m}\label{eq:Wiener_Integral}\end{alignat}
where $m$ is Wiener measure on $C_{0}^{t}$. Hence
\begin{equation}
\begin{split}
&\int_{0}^{t}e^{-(t-s)H_{0}}V(s)\Phi_{(R);l;l}^{1,s}(-H_{0},V)\, ds  =\\
&\int_{0}^{t}\sum_{m=0}^{\infty}a_{m}\int_{\Delta_{m}(s)}\int_{C_{0}^{t}}G(y+x)\, m(dy)dt_{1}\cdots dt_{m}ds  =\label{eq:RHS_Finish_Part_1}\\
&\int_{0}^{t}\int_{C_{0}^{t}}\left(\sum_{m=0}^{\infty}a_{m}\int_{\Delta_{m}(s)}G(y+x)\, dt_{1}\cdots dt_{m}\right)m(dy)ds
\end{split}
\end{equation}
where the last equality follows from the dominated convergence theorem
and Fubini's theorem. Finally, we may write the sum above in a different
way using the definition of $G$. In particular, all of the terms in
the product defining $G$ commute and we may write the integral of
$G$ over $\Delta_{m}(s)$ as
\begin{equation}
\frac{1}{m!}\left(\int_{0}^{s}V(y(u)+x,u)\, du\right)^{m}V(y(s)+x,s)\psi(y(t)+x).\label{eq:Unraveling_G}
\end{equation}
Hence, for Lebesgue almost every $x\in\mathbb{R}^{d}$, 
\begin{equation}
\begin{split}
&\int_{0}^{t}e^{-(t-s)H_{0}}V(s)\Phi_{(R);l;l}^{1,s}(-H_{0},V)\, ds=\label{eq:RHS_Finish_Part_2}\\
&\int_{0}^{t}\int_{C_{0}^{t}}\left(\sum_{m=0}^{\infty}\frac{a_{m}}{m!}\left(\int_{0}^{s}V(y(u)+x,u)\, du\right)^{m}\right)V(y(s)+x,s)\psi(y(t)+x)\, m(dy)ds 
\end{split}
\end{equation}
and our integral equation is, in this setting,
 \begin{multline}
f_{(R);l;l}^{t}(-H_{0},V)\psi(x)=g_{1}(0)e^{-tH_{0}}\psi(x)+\\
\int_{0}^{t}\int_{C_0^t}\left(\sum_{m=0}^{\infty}\frac{a_{m}}{m!}\left[\int_{0}^{s}V(y(u)+x,u)\, du\right]^{m}\right)V(y(s)+x,s)\psi(y(t)+x)\, m(dy)ds.\label{eq:Integral_Eqn_H0}
\end{multline}
If $V:\,\mathbb{R}^{d}\to\mathbb{R}$ is a time independent potential
with $V\in L^{\infty}(\mathbb{R}^{d})$, the previous calculations
go through in exactly the same was as above and we obtain
\begin{multline}
f_{(R);l;l}^{t}(-H_{0},V)\psi(x)=g_{1}(0)e^{-tH_{0}}\psi(x)+\\
\int_{0}^{t}\int_{C_0^t}\left(\sum_{m=0}^{\infty}\frac{a_{m}}{m!}\left[\int_{0}^{s}V(y(u)+x)\, du\right]^{m}\right)V(y(s)+x)\psi(y(t)+x)\, m(dy)ds.\label{eq:Time_Indep_IE}
\end{multline}
With this equation in hand, we are now in a position where we can investigate some examples
by making a specific choice for the  function $g_{1}$.

We will start in the time independent setting, taking $g_{1}(z_{1})=\frac{1}{1-z_{1}}$
and will assume here that $\| V\|_{\infty}<1$. As is well known, we may write
\[
g_{1}(z_{1})=\sum_{m=0}^{\infty}z_1^m.
\]
for $|z_1|<1$. We replace $V$ with $-V$ and obtain
\begin{multline}
f_{(R);l;l}^{t}(-H_{0},-V)\psi(x)=e^{-tH_{0}}\psi(x)-\\
\int_{0}^{t}\int_{C_{0}^{t}}\left(\sum_{m=0}^{\infty}\frac{(-1)^{m}}{m!}\left(\int_{0}^{s}V(y(u)+x)\, du\right)^{m}\right)
V(y(s)+x)\psi(y(t)+x)\, m(dy)ds\label{eq:Time_Indep_Spec_Case}
\end{multline}
where it is clear that $a_{m}=1$ and $g_{1}(0)=1$. Fubini's
theorem enables us to write this equation as 
\begin{equation}
\begin{split}
&f_{(R);l;l}^{t}(-H_{0},-V)\psi(x)=e^{-tH_{0}}\psi(x)-\\
&\int_{C_{0}^{t}}\int_{0}^{t}\exp\left[-\int_{0}^{s}V(y(u)+x)\, du\right]
 V(y(s)+x)\psi(y(t)+x)\, ds\, m(dy)\label{eq:IE_Integrated}\\
&=e^{-tH_{0}}\psi(x)-\int_{C_{0}^{t}}\left(1-\exp\left(-\int_{0}^{t}V(y(u)+x)\, du\right)\right)\psi(y(t)+x)\, m(dy)\\
&=e^{-tH_{0}}\psi(x)-\int_{C_{0}^{t}}\psi(y(t)+x)\, m(dy)+\int_{C_{0}^{t}}\exp\left(-\int_{0}^{t}V(y(u)+x)\, du\right)\cdot\\
&\psi(y(t)+x)\, m(dy)\\
&=e^{-tH_{0}}\psi(x)-e^{-tH_{0}}\psi(x)+e^{-t(H_{0}+V)}\psi(x)=e^{-t(H_{0}+V)}\psi(x)
\end{split}
\end{equation}
where we've applied the Feynman-Kac formula (see, for example, \cite{JLBook},
Chapter 12) to the third term after the third equality above and
the well known Wiener integral representation of the heat semigroup
to the second term after the third equality. Hence we see that the
reduced disentangling developed above supplies solutions to the heat
equation in this setting, for Lebesgue almost every $x\in\mathbb{R}^{d}$. 

We next make a connection between the reduced disentangling and various types of analytic Feynman integrals. We will stay in the setting of the calculation immediately
above, i.e. $g_{1}(z_{1})=\frac{1}{1-z_{1}}$ and $\| V\|_{\infty}<1$.
\begin{rem}
We briefly note that the example here involving $H_0$ and $V$ should also work with $H_0$ replaced by any operator $\alpha$ such that $-\alpha$ generates a $C_0$ semigroup and $V$ replaced by any bounded linear operator $A$ with uniform norm $\|A\|<1$. The proof would not appeal to Wiener measure.
\end{rem}
To begin, we
recall the definition of the analytic in time
Feynman integral $J^{t}(F)$ given on page 299 (Definition 13.2.1) of \cite{JLBook}:

\begin{defn}
Given $t>0$, $\psi\in L^{2}(\mathbb{R}^{d})$ and $\xi\in\mathbb{R}^{d}$,
consider the expression\[
(J^{t}(F)\psi)(\xi)=\int_{C_{0}^{t}}F(x+\xi)\psi(x(t)+\xi)\, m(dx).\]
The operator-valued function space integral $J^{t}(F)$ exists if
and only if the Wiener integral above defines $J^{t}(F)$ as an element
of $\mathcal{L}(L^{2}(\mathbb{R}^{d}))$. If $J^{t}(F)$ exists for
every $t>0$ and, in addition, has an extension (necessarily unique)
as a function of $t$ to an operator-valued analytic function on $\mathbb{C}_{+}$
(the set of complex numbers with positive real part) and a strongly
continuous function on $\bar{\mathbb{C}}_{+}$, we say that $J^{t}(F)$
exists for all $t\in\bar{\mathbb{C}}_{+}$. When $t$ is purely imaginary,
$J^{t}(F)$ is called \emph{the analytic in time operator-valued Feynman
integral}.
\end{defn}
In particular we will consider $J^{t}(F_{V})$ where $F_{V}$ is given
by
\begin{equation}
F_{V}(y)=\exp\left(-\int_{0}^{t}V(y(s))\, ds\right)\label{eq:Wiener_FNL}
\end{equation}
for $y\in C_{0}^{t}$. Indeed, as seen in Remark 13.2.2 on page 299
of \cite{JLBook},
\begin{equation}
e^{-t(H_{0}+V)}=J^{t}(F_{V}).\label{eq:Analytic_In_Time_FI}
\end{equation}
Hence we see that 
\begin{equation}\label{eq:Reduced_Heat_J}
\left[f_{(R);l;l}(-H_0,-V)\right]\psi(x)=J^t(F_V).
\end{equation}
Also, as is shown in Theorem 13.3.1 of \cite{JLBook},
\begin{equation}
J^{it}(F_{V})=e^{-it(H_{0}+V)}.\label{eq:Analytic_in_time_imaginary}
\end{equation} 
(So $J^t(F_V)$ is a solution to the heat equation for real $t$ and a solution to the Schrodinger equation for imaginary $t$.)
Moreover, if we replace $-H_{0}$ with $-iH_{0}$ and $-V$ with $-iV$
in the reduced disentangling $f_{(R);l;l}^{t}$ that we calculated
above, we obtain, with $\psi\in L^{2}(\mathbb{R}^{d})$, 
\begin{equation}
\left[f_{(R);l;l}^{t}(-iH_{0},-iV)\right]\psi(x)=e^{-it(H_{0}+V)}\psi(x)\label{eq:Imaginary_Reduced_Disent}\end{equation}
where the integrals appearing in the calculation must now be interpreted
in the mean. Hence, for $\psi\in L^{2}(\mathbb{R}^{d})$, 
\begin{equation}
\left[f_{(R);l;l}^{t}(-iH_{0},-iV)\right]\psi(x)=e^{-it(H_{0}+V)}\psi(x)=J^{it}(F_{V})\psi(x).\label{eq:Relation_Between}
\end{equation}
It should be noted that if $\psi\in L^{1}(\mathbb{R}^{d})\cap L^{2}(\mathbb{R}^{d})$,
then the integrals in the reduced disentangling as well as the integrals
in $J$ do not have to be interpreted in the mean. 
\begin{remark}
The class of potentials considered in Chapter 13 of \cite{JLBook} is quite large but does, of course, include the type of potential we are considering in this paper.
\end{remark}

We now move on to the setting where the potential $V$ is time dependent.
In this case we will make use of the mixed norm space $L_{\infty1;l}$.
We also need to recall the following definition (Definition 15.2.1
on page 410 of \cite{JLBook}) for the analytic (in mass) operator-valued
Feynman integral.

\begin{defn}
Fix $t>0$. Let $F$ be a function from $C_{0}^{t}$ to $\mathbb{C}$.
Given $\lambda>0$, $\psi\in L^{2}(\mathbb{R}^{d})$ and $\xi\in\mathbb{R}^{d}$,
we consider the expression\[
(K_{\lambda}^{t}(F)\psi)(\xi)=\int_{C_{0}^{t}}F(\lambda^{-1/2}x+\xi)\psi(\lambda^{-1/2}x(t)+\xi)\, m(dx).\]
The \emph{operator-valued function space integral $K_{\lambda}^{t}(F)$
exists for} $\lambda>0$ if the integral above defines $K_{\lambda}^{t}(F)$
as an element of $\mathcal{L}(L^{2}(\mathbb{R}^{d}))$. If, in addition,
$K_{\lambda}^{t}(F)$, as a function of $\lambda$, has an extension
(necessarily unique) to an analytic function on $\mathbb{C}_{+}$
and a strongly continuous function on $\bar{\mathbb{C}}_{+}\backslash\left\{ 0\right\} $,
we say that $K_{\lambda}^{t}(F)$ exists for $\lambda\in\bar{\mathbb{C}}_{+}\backslash\left\{ 0\right\} $.
When $\lambda$ is purely imaginary, $K_{\lambda}^{t}(F)$ is called
the \emph{analytic (in mass) operator-valued Feynman integral of}
$F$.
\end{defn}
For $y\in C_{0}^{t}$ and $V\in L_{\infty1;l}([0,T]\times\mathbb{R}^{d})$
we define (see equation \eqref{eq:Time_Indep_IE})
\begin{equation}
F_{V}(y)\,:=\sum_{m=0}^{\infty}\frac{(-1)^{m+1}a_{m}}{(m+1)!}\left(\int_{0}^{t}V(s,y(s))\, ds\right)^{m}.\label{eq:New_F}\end{equation}
That is, for \[
G(z)=\sum_{m=0}^{\infty}\frac{(-1)^{m+1}a_{m}}{(m+1)!}z^{m+1}\]
 with radius of convergence strictly greater than $\| V\|_{\infty1;l}$
we have \[
F_{V}(y)=G\left(\int_{0}^{t}V(s,y(s))\, ds\right)\]
and it follows that \[
\int_{C_{0}^{t}}F_{V}(y+x)\psi(y(t)+x)\, m(dy)\]
exists and is equal to $K_{1}^{t}(F_{V})\psi$. Recalling that $g_{1}(z_{1})=\sum_{m=0}^{\infty}a_{m}z_{1}^{m}$,
we can write \begin{equation}
\left[f_{(R);l;l}^{t}(-H_{0},-V)\right]\psi(x)=g_{1}(0)\left(e^{-tH_{0}}\right)\psi(x)-K_{1}^{t}(F_{V})\psi(x).\label{eq:Analytic_Op_Valued_Rel}\end{equation}

\begin{rem}
We have taken $g_{1}$ to be analytic in $|z_{1}|<\| V\|_{\infty1;l}$
and such that its series converges on $|z_{1}|=\| V\|_{\infty1;l}$.
Hence, for the series defining $G(z)$ the radius of convergence is
infinite, as $\lim_{m\to\infty}\left|\frac{a_{m}}{a_{m+1}}\right|=constant$ as the series for $g_1$ has a finite radius of convergence.  
\end{rem}

\end{document}